\documentclass[a4paper,12pt]{opelsarticle}
\usepackage{amsmath,amsfonts,graphicx,tikz,pgfplots, pdfsync}
\usepackage[colorlinks=true, citecolor=red]{hyperref}
\usepackage{fullpage, graphicx}
\usetikzlibrary{patterns}

\def\eq#1{\Blue{\begin{eqnarray}#1\end{eqnarray}}}

\DeclareGraphicsExtensions{.pdf, .jpg, .tif}

\definecolor{dbrown}{RGB}{0.41,0.05,0.05}
\definecolor{darkbrown}{RGB}{194,52,42}
\definecolor{darkred}{RGB}{250,30,30}
\def\Blue#1{{\color{blue}{#1}\color{black}}}
\def\Black#1{{\color{black}{#1}\color{black}}}

\def\E{{\mathbb E}}
\def\p{\partial}

\def\R{{\mathbb R}}

\def\L{{\mathcal L}}
\def\EE{{\mathcal E}}

\def\vone{{\mathbf 1}}

\def\vx{{\mathbf x}}
\def\vy{{\mathbf y}}
\def\vb{{\mathbf b}}
\def\vA{{\mathbf A}}

\def\ttdefault{cmtt}
\newtheorem{theorem}{Theorem}

\DeclareRobustCommand{\LaTeX}{L\kern-.26em%
        {\sbox\z@ T%
         \vbox to\ht\z@{\hbox{\check@mathfonts
           \fontsize\sf@size\z@
           \math@fontsfalse\selectfont
          A\,}%
         \vss}%
        }%
     \kern-.15em%
    \TeX}
\makeatother

\begin{document}

\title{Calibration of a Fluid-Structure Problem with Keras}
\author{
Olivier Pironneau \footnote{\emph{olivier.pironneau@gmail.com }, LJLL, Sorbonne University, Paris, France.}
}

\parindent=0pt
\begin{frontmatter}
\begin{abstract}
In this short paper we report on an inverse problem issued from a physical system, namely a fluid structure problem where the parameters are the rigidity constant, the solid-fluid density ratio and the fluid viscosity. We have chosen a simple geometry so that the numerical solution of the system is easy.  We compare the solution of this inverse problem by a  Neural Network with a more classical solution obtained with a genetic algorithm. The Neural Network does much better.
 \end{abstract}
\end{frontmatter}

\section{Introduction}
Let $\L_p w=f$ be the Partial Differential Equations (PDE) of a physical system which defines $w$ for a given set of parameters $p$; the inverse problem is to recover $p$ when $w$ is known.
Such inverse problems may be solved by optimization (see \cite{tarantola} ):
 \eq{
p= \hbox{argmin}_{p\in P_{ad}}\{ \|u^p-w\|^2~:~\L_p u^p=f\}
}
where $P_{ad}$ is the set of admissible values.  If we know only a finite number of values of $w$, say at $\{z_k\}_1^K$ then we may solve, with $w^k=w(z_k)$,
 \eq{\label{opt}
p= \hbox{argmin}_{p\in P_{ad}}\{\frac1K \sum_1^K|u^p(z_k)-w^k|^2~:~\L_p u^p=f\}
}
Deep Learning takes another approach; a neural network attempts to build an approximation of the inverse PDE and computes $\{w^k\}_1^K\mapsto p_{NN}$.  The neural network is defined by a set of parameters $\pi$. Training the network means adjusting the parameters $\pi$ so that $\{p_{NN}^m,\}_1^M$ is closest to  $\{p^m\}_1^M$ for $M$``samples'' $[\{w_m^k\}_1^K]_1^M$.  Once trained it is hoped that  the network  produce a good approximation $\bar p$ for a new set of data $\{\bar w^k\}_1^K$.
 
This kind of deep learning  \cite{goodfellow} is not just another inversion method because one has also the capacity to store $\pi$ and reuse the network in another context, another machine, like a tablet, and obtain inversion results very cheaply and fast.

In this paper we wish test the potential of  DNN (deep neural network)  to find the characteristics of an incompressible hyper-elastic material and of the fluid around it from the knowledge of a few data point in the space-time response to an external sollicitation.
The performance of the DNN inverter will be compared with a genetic algorithm inverter using CMA-ES\cite{hansen}.

\section{A simple Fluid-Strucuture Testing Device}

In 2 dimensions an incompressible hyper-elastic Mooney-Rivlin solid material (also called Neo-Hookian) is characterised by only one parameter, the Helmholtz potential constant $c_1$, besides its density $\rho^s$.  Guessing the viscosity of a viscous fluid is not much simpler; so there is an interest for a simple probe to measure all three.

Consider the following axisymmetric device:  a rigid fixed cylinder coated by an incompressible material , immersed in a viscous fluid (see figure \ref{disks}); the whole aparitus is inside an outer  cylinder.  The outer cylinder can rotate; it is an external  sollicitation  which induces motions in the fluid and the solid.    First the system is at rest and then a constant angular velocity is given to the outer cylinder.  This causes the fluid to rotate with an angular velocity  function of the radial distance $r$ to the main axis; in turn due to the friction of the fluid at the interface, the incompressible material will be have  an angular velocity $u_\theta$ which is also only a function of $r$  and time . Due to  elasticity $u_\theta$ will oscillate with time until numerical dissipation and fluid viscosity damps it. 

In a two dimensional cut perpendicular to the main axis, the velocities and displacements are two dimensional and axisymmetric as well. Hence
the geometry is a ring of inner and outer radii, $R_0$ and $R_1$, with the incompressible material between $R_0$ and $R$ and fluid between $R$ and $R_1$.  Because of the incompressibility of the fluid and axial symmetry, $R$ is constant.

\medskip

In this paper $R_0=3$, $R=4$, $R_1=5$.  The radial velocity of  the outer cylinder is also always kept at 3 rpm.   So the parameters of the problem are $c_1,\rho^s,\mu^f$.
\begin{figure}
\begin{minipage}[b]{0.45\linewidth}
\centering
\begin{tikzpicture}[xscale=0.9,yscale=0.9]
 \draw (3,0) ellipse (1.8 and 2.2) ;
 \filldraw[green] (3,0) ellipse (1.4 and 1.6) ;
 \draw[green,thick,densely dotted] (6,0) ellipse (1.4 and 1.6) ;
 \draw[green,thick,densely dotted](3,1.6)--(6,1.6);
 \draw[green,thick,densely dotted](3,-1.6)--(6,-1.6);
 \filldraw[white] (3,0) ellipse (0.9 and 1.1)  ;
 \filldraw[pattern=dots] (3,0) ellipse (0.9 and 1.1)  ;
 \draw[thick,densely dotted] (6,0) ellipse (0.9 and 1.1)  ;
  \draw[green,densely dotted](3,1.1)--(6,1.1);
  \draw[green,densely dotted](3,-1.1)--(6,-1.1);
  \draw(3,2.2)--(6,2.1);
 \draw(3,-2.2)--(6,-2.2);
 \draw (6,2.1) arc (90:-90:1.6 and 2.15);
 \draw (0.7,0.1) node [left] {{\tiny Fixed rod}}; 
 \draw (1.7,1.4) node [left] {{\tiny Hyperelastic material}};
 \draw[->,red,thick,densely dotted](1.7,1.4)--(2.8,1.4);
  \draw (0.8,1.7) node [left] {{\tiny Fluid}};
 \draw[->,red,thick,densely dotted](0.7,1.7)--(2.8,1.7);
 \draw (4.7, 1.8) node [right] {{\tiny outer cylinder}};
 \draw [->, red](5,0.5) arc (0:45: 2.1 and 1.8);
\draw [red,thick,densely dotted] (0.7,0)--(2.5,0);
\end{tikzpicture}
\caption{\label{disks}{\it  A fluid-structure system inside a rotating cylinder  with a fixed rod in its center .  Sketch of the system. }}
\end{minipage}
\hskip 1cm
\begin{minipage}[b]{0.45\linewidth}
\begin{tikzpicture}[scale=0.7]
\begin{axis}[legend style={at={(1,1)},anchor=north east}, compat=1.3,
  xmin=3., xmax=5,ymin=-0.2,ymax=5,
  xlabel= {Radial distance},
  ylabel= {Angular Velocity}]
\addplot[thick,solid,color=black,mark=none,mark size=1pt] table [x index=0, y index=1]{fsiplot.txt};
\addlegendentry{ $t=0.25$}
\addplot[thick,dashed,color=blue,mark=none,mark size=1pt] table [x index=0, y index=2]{fsiplot.txt};
\addlegendentry{ $t=1.25$}
\addplot[thick,dotted,color=red,mark=none,mark size=1pt] table [x index=0, y index=3]{fsiplot.txt};
\addlegendentry{ $t=2.5$}
\addplot[thick,dotted,color=black,mark=none,mark size=1pt] table [x index=0, y index=5]{fsiplot.txt};
\addlegendentry{ $t=5$}
\end{axis}
\end{tikzpicture}
\caption{\it A 2d calculation showing the angular velocities along the radial axis at several  times. }
\label{relerr} 
\end{minipage}
\end{figure}

As everything is axisymmetric the computation can be done  in polar coordinates $r,\theta$, and the fluid-solid system reduces to
\eq{\label{rotdisk}
\rho\p_t u_\theta - \frac1r\p_r[\xi^f r\p_r u_\theta+\xi^s r\p_r d] =0,~~
\p_t d = u_\theta,~r\in(R_0,R_1),~~v_{|R_0}=0,~{u_\theta}_{|R_1}=3,
}
with $\rho=\rho^s\vone_{r\leq R} + \rho^f\vone_{r>R}$, $\xi^s=2c_1\vone_{r\leq R}$, $\xi^f= {\mu}\vone_{r>R}$, and with $d(r,0)=0$.

The system is approximated in time by an implicit Euler scheme and in space by the finite element method of degree 2.  The linear systems are solved by LU factorization.  All parameters are shown in Appendix \ref{appdxA} which is a freefem \cite{freefem} implementation.  Figure \ref{relerr} shows results at 4 time instants.

Comparison between this one dimensional system and the full Eulerian 2D fluid-structure system has been reported in \cite{FHOP}

\section{Solution with a Genetic Algorithm}
Consider (\ref{opt}) with $K=10\times 5$ data points and $w^k=u_\theta(r^i,t^j)$ for 10 (resp 5) uniformly distributed $r^i$ (resp.  $t^j$) , CMAES \cite{hansen} is a global stochastic optimizer  which, in our case, stops when it  reaches the ``tolerance'' $10^{-6}$ for the function evaluations. In the process it has done 569 function evaluations and the relative error on the approximate solution is
$[0.087\%,0.17\%,0.085\%]$ for $[c_1,\rho^s,\mu^f]$.
With $5\times5$ data points, the relative error is $[0.13\%,0.08\%,0.15\%]$. 

 It is computationally expensive (15min on a Mac book pro) but it indicates that the problem has a solution and the precision can be improved to any amount, it seems.

\section{Solution with a Neural Network}

A neural network is characterized by the number of hidden layers (depth) and the number of neurons in each layer (width).  A neuron takes an input vector $\vx$ and produces an output vector $\vy$ with $y_i = \phi([\vA\vx+\vb]_i)$; $\vA$ is a rectangular matrix, $\vb$ is a constant vector and $\phi$ is a non-linear function the so-called ``activation''.  Here we have used systematically \texttt{ReLU}$(x):=x^+$ (rectified linear unit).

The DNN-parameters (called $\pi$ in the introduction),  $\vA,\vb$ of each hidden layer are adjusted by a stochastic gradient method to minimize a ``loss'' function. 

The loss function is the least square error between the output of the DNN and the output of the PDE.  Hence $M$ solutions of the PDE are computed, each with its own parameters $p=(c_1,\rho^s,\mu^f)$. For the training phase
the data of each sample consists of $p$ and $K$ values of the solution of the PDE.

The down side of deep  learning is that a large number of samples required and for each sample the PDE needs to be solved.  A few hours are needed to compute 5000 samples while the learning phase takes a few tens of seconds and the test phase microseconds.

A Keras code is given in Appendix \ref{appdxB}.  Note the simplicity of the \texttt{Python}  program.

\subsection{Influence of the number of Samples}
We have taken a DNN of width 100 and depth 1.  Training is done with $M$ samples each containg $K=10\times 5$ time-space values .  The influence of the number of samples $M$ on the precision of the recovered parameters $c_1,\rho^s$ and $\mu^f$ is shown on Table \ref{tab5} and Figure \ref{relerr}.  Table \ref{tab6} illustrates the precision with some values given by the DNN and compared with the true values.

\begin{table}[htp]
\begin{center}
\caption{\textbf{Precision} versus the number of samples:  average absolute (left) and relative (in \%)  errors.\label{tab5}
}
{\footnotesize\begin{tabular}{|c||c|c|c||c|c|c|}
\hline
Samples& $\|c_1-{c_1}_{true}\|$&$\|\rho^s-{\rho^s}_{true}\|$&$\|\mu^f-{\mu^f}_{true}\|$& error$(c_1)/c_1$&error$(\rho^s)/\rho^s$&error$(\mu^f)/\mu^f$\cr
\hline 
250 &   0.1110471 & 0.1489467 & 0.06382506&2.8\% &3.8\% &1.7\% \cr 
500 & 0.0506839&  0.08803635 &0.03628453&1.2\%& 2.2\%& 0.9\% \cr 
1000& 0.04044557 &0.0534637  &0.02502107&1.1\%&1.4\%& 0.70\%\cr 
2000 & 0.02899077& 0.04001569& 0.01182538&0.89\%& 0.97\% &0.36\% \cr 
\hline
\end{tabular}}
\end{center}
\end{table}%

\begin{table}[htp]
\caption{\textbf{Example of results with a DNN} with 100 neurons and using 1000 samples }
\begin{center}
{\footnotesize\begin{tabular}{|c|c|c||c|c|c|}
\hline
${c_1}_{NN}$&${\rho^s}_{NN}$&${\mu^f}_{NN}$&${c_1}_{true}$&${\rho^s}_{true}$&${\mu^f}_{true}$\cr
\hline 
1.4074545 &1.0255045& 1.1367583 & 1.42379 &1.04122& 1.15782 \cr
0.64601064 &0.93555903 &1.4175779  & 0.620108& 0.905501& 1.42252  \cr
0.6114049 &1.1150271 &1.44915   & 0.611026 &1.10804 & 1.46802  \cr
[0.7821889 & 1.2131171 & 0.49625444 & 0.751256 &1.14957 & 0.505141 \cr
[0.69889516& 1.10933  &  1.3671701  & 0.705757 &1.13354 & 1.37266 \cr
\hline
\end{tabular}}
\end{center}
\label{tab6}
\end{table}%

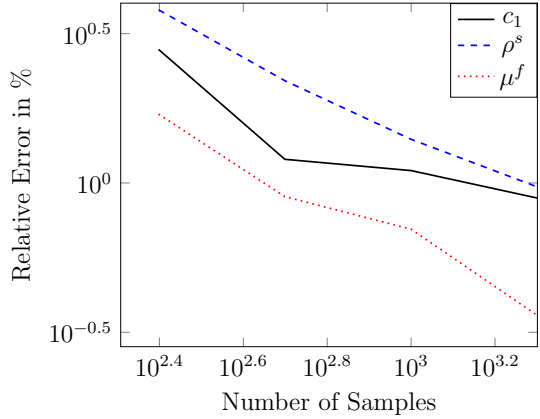
\begin{figure}[h!]
\centering
\begin{tikzpicture}[scale=0.8]
\begin{loglogaxis}[legend style={at={(1,1)},anchor=north east}, compat=1.3,
  xmin=0., xmax=2000,ymin=0,ymax=4,
  xlabel= {Number of Samples},
  ylabel= {Relative Error in \%}]
\addplot[thick,solid,color=black,mark=none,mark size=1pt] table [x index=0, y index=1]{precisionfsidnn.txt};
\addlegendentry{ $c_1$}
\addplot[thick,dashed,color=blue,mark=none,mark size=1pt] table [x index=0, y index=2]{precisionfsidnn.txt};
\addlegendentry{ $\rho^s$}
\addplot[thick,dotted,color=red,mark=none,mark size=1pt] table [x index=0, y index=3]{precisionfsidnn.txt};
\addlegendentry{ $\mu^f$}
\end{loglogaxis}
\end{tikzpicture}
\caption{Relative error versus the number of samples (log-log scales).}
\label{relerr} 
\end{figure}

\subsection{Influence of the Number of Data Points}

With 5 points in $r$-axis and $2$ points in time, the relative precision is $[2.4\%, 5.2\%, 1.5\%]$.
With 10 points in $r$-axis and $2$ points in time, the relative precision is $[2.0\%, 6.0\% 1.3\%]$.
With 5 points in $r$-axis and $5$ points in time, the relative precision is $[1.1\%, 2.1\%, 0.8\%]$.
With 10 points in $r$-axis and $5$ points in time, the relative precision is $[1.1\%, 0.97\%, 0.36\%]$.

In all 4 computations the number of samples is 1000.

\subsection{Influence of the Network Architecture}
We consider the case $10\times 5$ data points and 1000 samples.
\begin{itemize}
\item With 1 layer of width 1000, the number of epochs to trigger \emph{early stopping} is 111 and the loss is $0.166689$.
\item With 2 layers of width 100, the number of epochs  is 196 and the loss is $0.10234474$.
\item With 1 layer of width 100, the number of epochs is 374 and the loss is $0.13168786$.
\item With a 2 layer of width 50, the number of epochs is 596 and the loss is $0.06836466$.
\item With 3 layers of width 34, the number of epochs is 393 and the loss is  $0.10246510$.
\item With a 2 layer of width 50 and a layer of width 10 in between, the number of epochs is 888 and the loss is $0.06401216$.
\end{itemize}
In view of the fact that all these numbers are  realisations of random processes, it is not easy to draw any conclusion, but it seems to indicate that the architecture of the network is not very important, despite the following results:

\begin{theorem} \cite{yarotsky}
Functions of $W^{n,\infty}([0,1]^d)$ can be $\delta$-approximated by ReLU networks with depth $O(\log\frac1\delta)$ and width $O(\delta^{-\frac d n})$.
\end{theorem}

\begin{figure}[htbp]
\begin{center}
\includegraphics[width=10cm]{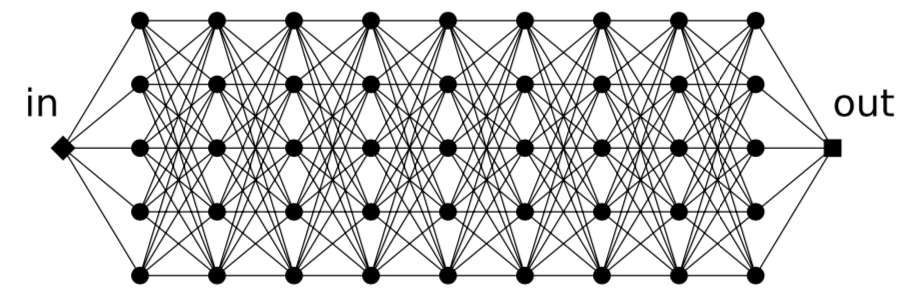}
\caption{A Deep network with constant width and large depth.}
\label{figfive}
\end{center}
\end{figure}

\begin{theorem}\cite{yarotsky2}  {A bounded function of $x\in \R$ with bounded first derivative can be approximated with precision $O(1/(N\log N))$ by the network of Figure \ref{figfive} of fixed width (here 5) and of depth N.}
 \end{theorem}
 
 \subsection{Critical Convergence}
 Let us explain why the precision is hard to adjust. Assuming  $M$ samples
 $\{x_i,u^\epsilon(x_i)\}_1^M$  known with stochastic error $\epsilon$ to simulate 
\eq{
\R^d\ni x \mapsto u(x)\in\R^s \hbox{\Black{ by a DNN }}x\mapsto U(x) 
}
Then the most sensitive thing to do is to 
\begin{itemize}
\item Rescale $x_i\in(-1,1)^d$ and to 
\item Minimize with respect to the parameters which define $U$ the loss function  using a DNN of \Blue{width and depth $M_d$}, and observe the loss
 \eq{
 \EE_N(U)=\frac1M\sum_{i=1}^M \| U({x}_i)-u^\epsilon({x}_i)\|^2
 }
\item The AdaGrad-Norm algorithm \cite{adagrad} will converge at rate $O(\frac{log P}P)$, after $P$ iterations,  to a local minimum $U^*$.
 \end{itemize}
 Then, for some $s>1$,
 \eq{
 \|\E[u^\epsilon|x]-U^*(x)\|\leq C M_d^{-s} + \EE_N(U).
 }
Hence the value of the ``loss" (the last term above) at convergence will control the precision because the first term on the right can be decreased to any small number by increasing the number of neurons.  Unfortunately the magnitude fo this last term can only be discovered a posteriori.

 \section{Conclusion}
The physical system chosen here is sufficiently simple to be sure of the accuracy of the numerical method. This has allowed us to concentrated on the parameters of the DNN.  The bottleneck is the number of samples needed to train the network, although here a few hundreds are sufficient for a precision around 1\%.  Consequently the DNN performs very well.  The architecture of the network is not important; the size of the data in each sample should also not be too small.  

We intend to build on this experience to solve a more difficult FSI problem with a free boundary as in \cite{glow80}.

\newpage

\section*{Appendix: Programs}
\subsection{ freefem program to generate the data} \label{appdxA} 
The freefem PDE solver for Mac, PC, Linux, is available for free download at texttt{www.freefem.org}. 

\includegraphics[width=13cm]{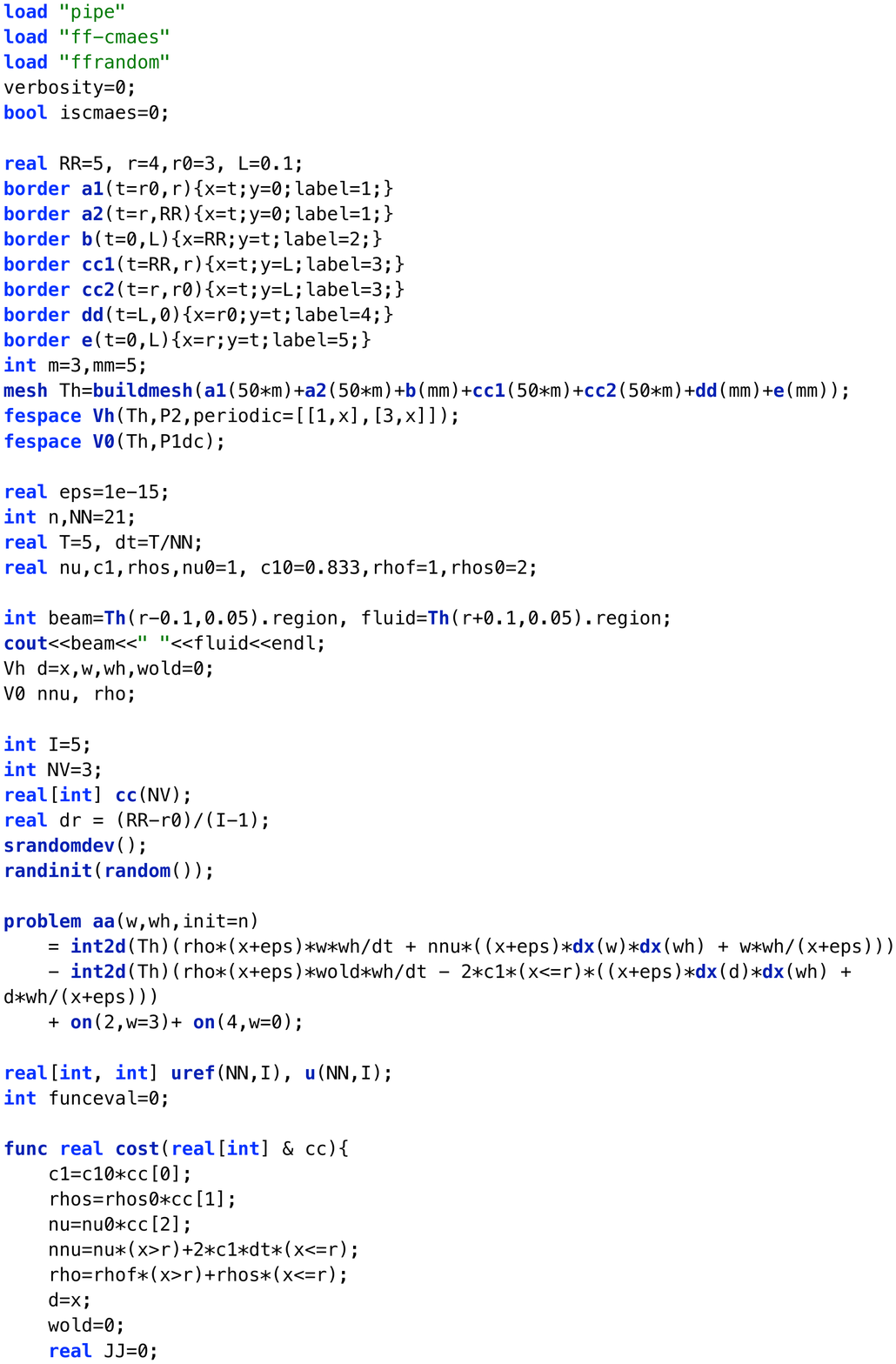}
~
\vskip-1.5cm
\includegraphics[width=13cm]{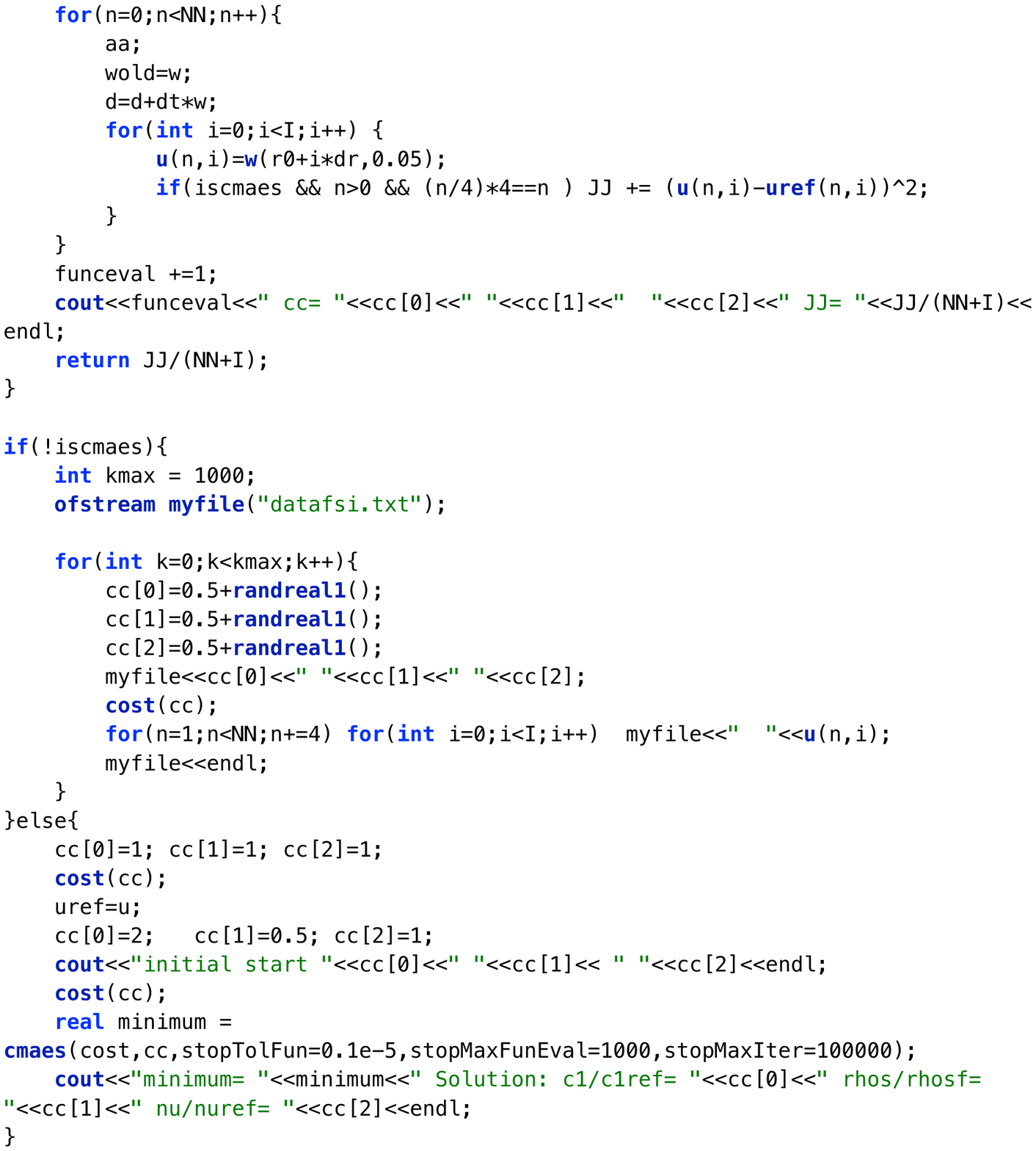}
\vskip-1.5cm

\subsection{ The Keras/Python program to compute $c_1$,$\rho^s$ and $\mu^f$ with a neural network}\label{appdxB}
\vskip-1.5cm
\includegraphics[width=13cm]{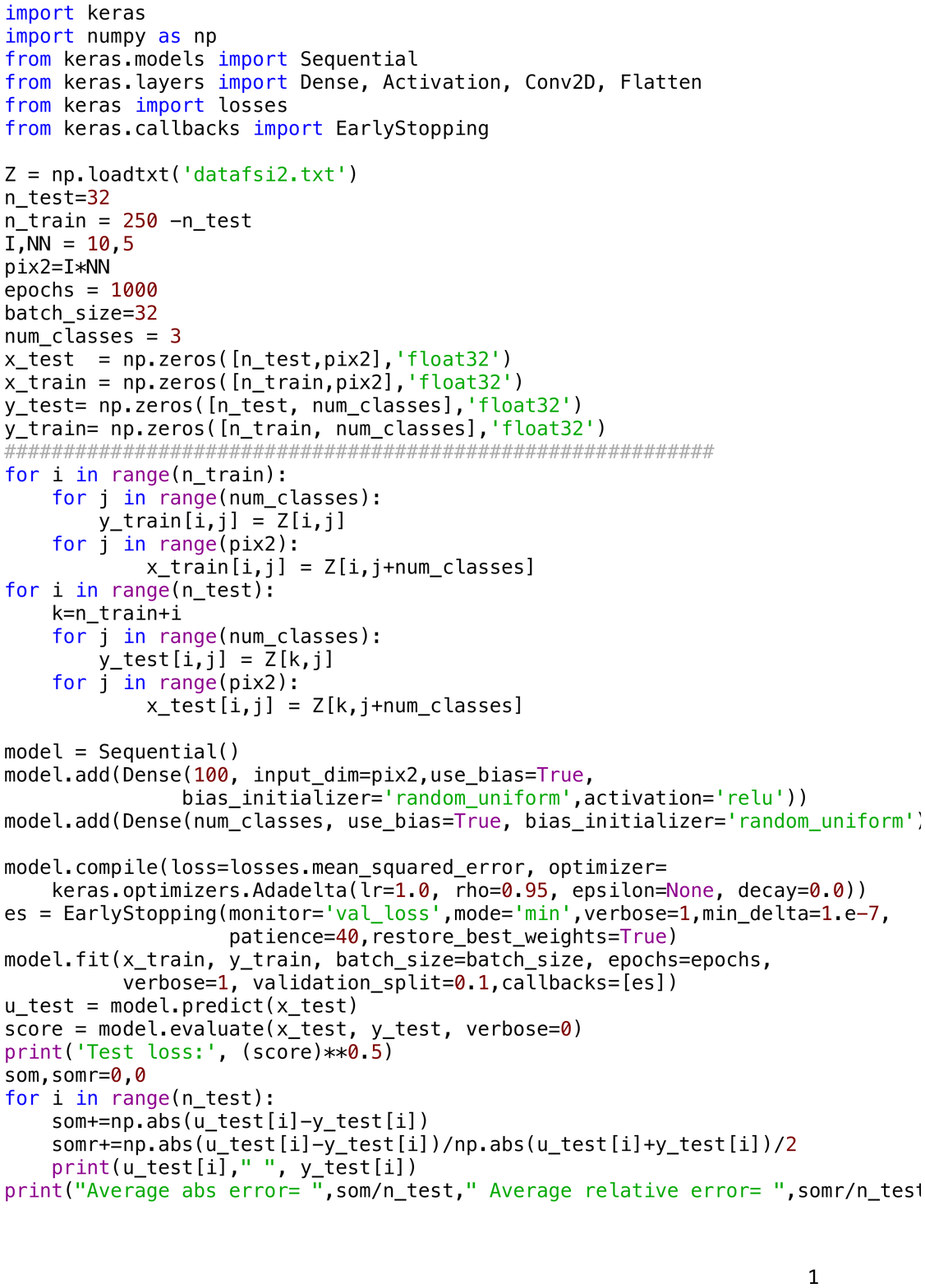}

\end{document}